# On a six-parameter generalized Burr XII distribution

A.K. Olapade

Department of Mathematics, Obafemi Awolowo University, Ile-Ife, Nigeria.

E-mail: akolapad@oauife.edu.ng

**Abstract**

In this paper, we derive a probability density function that generalizes the Burr XII distribution. The cumulative distribution function and the $n^{th}$ moment of the generalized distribution are obtained while the distribution of some order statistics of the distribution are established. A theorem that relate the new distribution to another statistical distribution is established.

## 1 Introduction.

A continuous random variable $X$ is said to follow a Burr XII distribution if its cumulative distribution function is given by

$$F(x) = 1 - (1 + \theta x^p)^{-m}, \quad 0 \leq x < \infty, \; m > 0, \; \theta > 0, \; p > 0; \tag{1.1}$$

and its probability density function is given by

$$f(x) = \frac{m\theta p x^{p-1}}{(1+\theta x^p)^{m+1}}, \quad 0 \leq x < \infty, \; m > 0, \; \theta > 0, \; p > 0. \tag{1.2}$$

If the location parameter $\mu$ and the scale parameter $\sigma$ are introduced in the equation (1.2), we have

$$f(x;\mu,\sigma,m,\theta,p) = \frac{m\theta p(\frac{x-\mu}{\sigma})^{p-1}}{\sigma[1+\theta(\frac{x-\mu}{\sigma})^p]^{m+1}}, \quad 0 \leq \mu \leq \infty, \; \sigma > 0. \tag{1.3}$$

Hence equation (1.3) is a five-parameter generalized Burr XII distribution.

Many reasearch work have been done on this distribution. Tadikamalla (1980) studied this distribution and other Burr type distributions. Sarah and Pushkarna (1999) obtained moments of order statistics from doubly truncated Lomax distribution which is a special case of Burr XII distribution. Begum and Parvin (2002) obtained moments of order statistics from doubly truncated Burr XII distribution while Haseeb and Khan (2002) listed Burr XII among a general class of distributions that satisfies $F(x) = ah(x) + b$ which they characterized by cosidering conditional moments of function of order statistics.

The Burr XII distribution is an important distribution because it has many other distributions like Pareto II or Lomax distribution (see Arnold (1983), Balakrishnan and Ahsanullah (1994)), log-logistic distribution, compound-Weibull or Weibull-Gamma and Weibull Exponential as particular cases of this distribution. Therefore, expressions for these distributions could be obtained from expressions of Burr XII distribution [Pareto or Lomax at $p = 1$ (Patil and Taillie (1994)), log-logistic at $m = 1$, $\theta = 1$ (Balakrishnan, Malik and Puthenpura (1987)), Weibull-Gamma at $\theta = 1/\delta$ (Tadikamalla (1980))]. So generalizing Burr XII distribution is generalizing all other distributions which come under special cases of Burr XII distribution.

## 1.1  Motivation for generalizing Burr XII distribution:

As mentioned in Ojo and Olapade (2005), It is well known in general that a generalized model is more flexible than the ordinary model and it is preferred by many data analysts in analyzing statistical data. Moreover, it presents beautiful mathematical exercises and broadened the scope of the concerned model being generalized. We like to mention that through generalizations, the ordinary logistic distribution with only two parameters (location and scale) has been developed into type I, type II and type III generalized logistic distributions which has three parameters each as shown in Balakrishnan and Leung (1988). Also, George and Ojo (1980) have through generalization developed a four-parameter generalized logistic distribution. Wu et-al (2000) have obtained a five-parameter generalized logistic distribution while Olapade (2004) developed an extended type I generalized logistic distribution which contains four parameters and generalizes the type I generalized logistic distribution of Balakrishnan and Leung (1988). Furthermore, Olapade (2005) presented a negatively skewed generalized logistic distribution which contains four parameters and further generalizes the type II generalized logistic distribution of Balakrishnan and Leung (1988). A six-parameter generalized logistic distribution was presented by Ojo and Olapade (2004). All these generalizations have become useful tools in the hand of researchers. In the same manner, we are of the opinion that generalizing Burr XII distribution will yield a similar result.

In this paper, a step forward is taken by defining a suitable random variable, the probability density function of a six-parameter generalized Burr XII distribution is derived. The cumulative distribution function is obtained and its $n^{th}$ moment is established. The probability density function of some of its order statistics are obtained.

## 2 The six-parameter generalized Burr XII distribution, its cumulative distribution function and moments

### 2.1 The six-parameter generalized Burr XII distribution

**Theorem 2.1.** *Suppose $Y_1$ and $Y_2$ are independently distributed continuous random variables. If $Y_1$ has an exponential density function*

$$f(y_1; \theta) = \theta e^{-\theta y_1}, \quad y_1 > 0, \quad \theta > 0$$

*and $Y_2$ has a gamma distribution with probability density function*

$$f(Y_2; m, \lambda) = \frac{\lambda^m}{\Gamma m} y_2^{m-1} e^{-\lambda y_2}, \quad y_2 > 0, \ m > 0, \ \lambda > 0. \tag{2.1}$$

*Then, the random variable*

$$X = \sqrt[p]{\frac{Y_1}{Y_2}}$$

*has the six-parameter generalized Burr XII distribution with parameters ($\mu = 0, \sigma = 1, \lambda, \theta, m, p$).*

**Proof.** The joint probability density function of $Y_1$ and $Y_2$ is

$$f(y_1, y_2) = \frac{\lambda^m \theta}{\Gamma m} y_2^{m-1} e^{-(\theta y_1 + \lambda y_2)}, \quad y_1 > 0, \ y_2 > 0, \ m > 0, \ \lambda > 0, \theta > 0. \tag{2.2}$$

Let $x_1 = \sqrt[p]{\frac{y_1}{y_2}}$ and $x_2 = y_2$. We obtain the probability density of the random variable $X_1$ as

$$f(x_1) = \frac{\lambda^m \theta p}{\Gamma m} \int_0^\infty x_1^{p-1} x_2^m e^{-x_2(\theta x_1^p + \lambda)} dx_2$$

$$= \frac{\lambda^m m \theta p x_1^{p-1}}{(\lambda + \theta x_1^p)^{m+1}}, \quad x_1 > 0, \ \lambda > 0, \ \theta > 0, \ m > 0, \ p > 0. \tag{2.3}$$

If we introduce the location parameter $\mu$ and the scaled parameter $\sigma$ in the equation (2.3) we have

$$f(x;\mu,\sigma,\lambda,\theta,m,p) = \frac{\lambda^m m\theta p(\frac{x-\mu}{\sigma})^{p-1}}{\sigma[\lambda + \theta(\frac{x-\mu}{\sigma})^p]^{m+1}}, \quad x > 0,\ \mu > 0,\ \sigma > 0,\ \lambda > 0, \theta > 0,\ m > 0,\ p > 0. \tag{2.4}$$

This probability density function in equation (2.4) is what we refer to as a six-parameter generalized Burr XII distribution. For the rest of this paper, we shall assume that $\mu = 0$ and $\sigma = 1$ without loss of generality.

## 2.2 Cumulative distribution function (cdf) of the generalized Burr XII distribution.

If a random variable $X$ has the generalized Burr XII probability density function given in the equation (2.3), the cdf of $X$ is given as

$$F(x) = \lambda^m m\theta p \int_0^x \frac{t^{p-1}}{(\lambda + \theta t^p)^{m+1}} dt$$

$$= 1 - \lambda^m(\lambda + \theta x^p)^{-m}, \quad \lambda > 0,\ \theta > 0,\ m > 0,\ p > 0. \tag{2.5}$$

The probability that a generalized Burr XII random variable $X$ lies in the interval $(\alpha_1, \alpha_2)$ is given as

$$Pr(\alpha_1 < X < \alpha_2) = \lambda^m[(\lambda + \theta\alpha_1^p)^{-m} - (\lambda + \theta\alpha_2^p)^{-m}], \text{ for } \alpha_1 < \alpha_2. \tag{2.6}$$

Hence given the value of the parameters $\lambda, \theta, p, m$ and an interval $(\alpha_1, \alpha_2)$, the probability $Pr(\alpha_1 < X < \alpha_2)$ can be easily computed.

## 3 Moments of the generalized Burr XII distribution

The $n^{th}$ moment of a random variable $X$ that has the generalized Burr XII distribution is

$$E[X^n] = \int_{-\infty}^{\infty} x^n f(x) dx \tag{3.1}$$

$$= \lambda^m m\theta p \int_0^\infty \frac{x^{n+p-1}}{(\lambda + \theta x^p)^{m+1}} dx = \frac{\theta m p}{\lambda} \int_0^\infty \frac{x^{n+p-1}}{(1 + \frac{\theta}{\lambda}x^p)^{m+1}} dx. \tag{3.2}$$

Let $z = \theta x^p/\lambda$, then $x = \sqrt[p]{\frac{\lambda}{\theta}z}$. So

$$E[X^n] = \frac{\lambda^{n/p}m}{\theta^{n/p}} \int_0^\infty \frac{z^{n/p}}{(1+z)^{m+1}} dz = \frac{\lambda^{n/p}m}{\theta^{n/p}} B(n/p + 1, m - n/p) \tag{3.3}$$

$$= m(\frac{\lambda}{\theta})^{n/p}\frac{\Gamma(n/p+1)\Gamma(m-n/p)}{\Gamma(m+1)} = (\frac{\lambda}{\theta})^{n/p}\frac{\Gamma(n/p+1)\Gamma(m-n/p)}{\Gamma m}, \quad m > n/p. \quad (3.4)$$

When $n = 1$ in equation (3.4), we obtain the mean of the generalized Burr XII distribution as

$$E[X] = (\frac{\lambda}{\theta})^{1/p}\frac{\Gamma(1/p+1)\Gamma(m-1/p)}{\Gamma m}, \quad m > 1/p. \quad (3.5)$$

Also, when $n = 2$,

$$E[X^2] = (\frac{\lambda}{\theta})^{2/p}\frac{\Gamma(2/p+1)\Gamma(m-2/p)}{\Gamma m}, \quad m > 2/p. \quad (3.6)$$

Hence, the variance of the generalized Burr XII distribution is obtained as

$$\sigma_X^2 = (\frac{\lambda}{\theta})^{2/p}[\frac{\Gamma(2/p+1)\Gamma(m-2/p)}{\Gamma m} - (\frac{\Gamma(1/p+1)\Gamma(m-1/p)}{\Gamma m})^2]. \quad (3.7)$$

# 4 Order Statistics of the generalized Burr XII distribution.

Let $X_1, X_2, ..., X_n$ be $n$ independently continuous random variables from the generalized Burr XII distribution and let $X_{1:n} \leq X_{2:n} \leq ... \leq X_{n:n}$ be the corresponding order statistics. Let $F_{X_{(r:n)}}(x)$, $(r = 1, 2, ..., n)$ be the cumulative distribution function of the $r^{th}$ order statistics $X_{(r:n)}$ and $f_{X_{(r:n)}}(x)$ denotes its probability density function. David (1970) gives the probability density function of $X_{(r:n)}$ as

$$f_{X_{(r:n)}}(x) = \frac{1}{B(r, n-r+1)} P^{r-1}(x)[1 - P(x)]^{n-r} p(x). \quad (4.1)$$

For the generalized Burr XII distribution with probability density function and cumulative distribution function given in the equations (2.3) and (2.5) respectively, by substituting $f(x)$ for $p(x)$ and $F(x)$ for $P(x)$ in the equation (4.1), we have

$$f_{X_{(r:n)}}(x) = \frac{1}{B(r, n-r+1)}[1 - \frac{\lambda^m}{(\lambda+\theta x^p)^m}]^{r-1}[\frac{\lambda^m}{(\lambda+\theta x^p)^m}]^{n-r}\frac{\lambda^m m\theta p x^{p-1}}{(\lambda+\theta x^p)^{m+1}}$$

$$= \frac{\lambda^{m(n-r+1)}m\theta p}{B(r, n-r+1)}\frac{x^{p-1}[(\lambda+\theta x^p)^m - \lambda^m]^{r-1}}{(\lambda+\theta x^p)^{mn+1}}. \quad (4.2)$$

## 4.1 The probability density function of the minimum and maximum observations from the generalized Burr XII distribution.

The minimum observation is denoted as $X_{1:n}$ and its probability density function is obtained by making $r = 1$ in the equation (4.2) to have

$$f_{X_{(1:n)}}(x) = \frac{\lambda^{mn} mn\theta p x^{p-1}}{(\lambda + \theta x^p)^{mn+1}}, \quad x > 0, \lambda > 0, \theta > 0, m > 0; \tag{4.3}$$

which is another generalized Burr XII distribution with parameter $(\lambda, \theta, mn, p)$.

The maximum observation is denoted by $X_{n:n}$ and its probability density function is obtained by making $r = n$ in the equation (4.2) to have

$$f_{X_{(n:n)}}(x) = \frac{\lambda^m mn\theta p x^{p-1}\{(\lambda + \theta x^p)^m - \lambda^m\}^{n-1}}{(\lambda + \theta x^p)^{mn+1}}, \quad x > 0, \lambda > 0, \theta > 0, p > 0, \ m > 0. \tag{4.4}$$

The $q^{th}$ moment of the minimum observation from the generalized Burr XII distribution is

$$E[X_{(1:n)}^q] = \int x^q f_{X_{(1:n)}}(x) dx$$

$$= \lambda^{mn} mn\theta p \int_0^\infty \frac{x^{p+q-1}}{(\lambda + \theta x^p)^{mn+1}} dx = mn(\frac{\lambda}{\theta})^{q/p} B(\frac{q}{p} + 1, mn - \frac{q}{p}), \tag{4.5}$$

where $B(.,.)$ is a complete beta function.

Hence, the mean of the minimum observation from the generalized Burr XII distribution is

$$E[X_{(1:n)}] = mn(\frac{\lambda}{\theta})^{1/p} B(\frac{1}{p} + 1, mn - \frac{1}{p}), \tag{4.6}$$

while the variance is

$$Var[X_{(1:n)}] = mn(\frac{\lambda}{\theta})^{2/p}[B(\frac{2}{p} + 1, mn - \frac{2}{p}) - mn(B(\frac{1}{p} + 1, mn - \frac{1}{p}))^2]. \tag{4.7}$$

# 5 A relationship between the generalized Burr XII and other statistical distribution.

In this section, we shall state and prove a theorem that relates the generalized Burr XII distribution to another statistical distribution.

**Theorem 5.1:** *Suppose $Y$ is a continuously distributed random variable with probability density function $f_Y(y)$, then the random variable*

$$X = \sqrt[p]{e^{\frac{Y}{m}} - \frac{\lambda}{\theta}}$$

*has a generalized Burr XII distribution with parameters $(p, m, \theta, \lambda)$ if $Y$ is exponentially distributed.*

**Proof:** The probability density function of an exponential random variable $Y$ is

$$f_Y(y) = e^{-y}, \quad y > 0.$$

By omitting all constants, the density of $X$ can be written as

$$f_X(x) \propto \frac{x^{p-1}}{(\lambda + \theta x^p)^{m+1}}. \tag{5.1}$$

Since any density function proportional to the right hand side of the equation (5.1) is that of a generalized Burr XII random variable, the proof is complete.

# References

[1] B.C. Arnold. *The Pareto Distribution.* International Co-operative Publishing House, Fairland, MD. (1983)

[2] N. Balakrishnan and M. Ahsanullah. *Relations for single and product moments of record values from Lomax distrbution.* Sankhya: The Indian Journal of Statistics. Vol. 56, Series B, Pt. 2, (1994), 140-146.

[3] N. Balakrishnan and M. Y. Leung, *Order statistics from the Type I generalized Logistic Distribution*, Communications in Statistics - Simulation and Computation. Vol. 17(1), (1988), 25-50.

[4] N. Balakrishnan, H.J. Malik and S. Puthenpura, *Best linear unbiased estimation of location and scale parameters of the log-logistic distribution.* Commn. Statist. - theory meth., 16(12), (1987), 3477-3495.

[5] A.A. Begum and S. Parvin, *Moments of Order Statistics from Doubly Truncated Burr Distribution*, J. Statist. Reasearch, 36(2), (2002), 179-190.

[6] H.A. David, *Order Statistics.* John Wiley, New York (1970).

[7] Athar Haseeb and A.H. Khan, *On Characterization of Probability Distributions through Conditional Expectation of Order Statistics*, J. Statist. Reasearch, 36(2), (2002), 131-136.

[8] E.O. George and M.O. Ojo, *On a generalization of the logistic distribution*, Annals of Statistical Mathematics, 32, 2, A, (1980), 161-169.

[9] M.O. Ojo and A.K. Olapade, *On a Six-parameter Generalized Logistic Distributions*, Kragujevac J. Math. Serbia. 26, (2004), 31-38.

[10] M.O. Ojo and A.K. Olapade, *On a generalization of the Pareto distribution*, Proceeding of the International Conference in honor of Prof. E.O. Oshobi and Dr. J.O. Amao, (2005), 65-70.

[11] A.K. Olapade, *On Extended Type I generalized logistic distribution*, International Journal of Mathematics and Mathematical Sciences, 57, (2004), 3069-3074.

[12] A.K. Olapade *On Negatively Skewed Extended Generalized Logistic Distribution*, Kragujevac J. Math. Serbia. 27, (2005), 175-182.

[13] G.P. Patil and C. Taillie, *Statistical evaluation of the attainment of interim cleanup standards at hazardous waste sites*. Environmental and Ecological Statistics **1**, (1994), pp. 333-351.

[14] J. Saran and N. Pushkarna, *Moments of order statistics from doubly truncated Lomax distribution*. J. Statist. Reasearch, 33(1), (1999), 57-66.

[15] P.R. Tadikamalla, *A look at the Burr and related distributions*. International Statist. Rev. 48, (1980), 337-344.

[16] Wu Jong-Wuu, Hung Wen-Liang and Lee Hsiu-Mei, *Some Moments and Limit Behaviors of the Generalized Logistic Distribution with Applications*, Proc. Natl. Sci. Counc. ROC(4) Vol. 24, No. 1, (2000), 7-14.